\documentclass[twocolumn]{autart}
\usepackage{makeidx}  
\usepackage{graphicx} 
\usepackage{mathrsfs}
\usepackage{arydshln}
\usepackage[tbtags]{amsmath}         
\usepackage{amssymb}                 
\usepackage{epsfig}                  
\usepackage{latexsym}                
\usepackage{array}                  
\usepackage{ifthen}                 
\usepackage{epstopdf}
\usepackage{subfigure}
\usepackage{theorem}
\usepackage{color}
\usepackage{comment}
\usepackage{CJK}
\usepackage{epstopdf}
\usepackage{float}
\newtheorem{theorem}{Theorem}
\newtheorem{lemma}{Lemma}
\newtheorem{assumption}{Assumption}
\newtheorem{remark}{Remark}

\newcommand{\bproof}{\noindent{\it Proof: }}
\newcommand{\eproof}{\hfill\rule{2mm}{2mm}}

\newcommand{\R}{\mbox{$\mathbb{R}$}}
\newcommand{\C}{\mbox{$\mathbb{C}$}}

\renewcommand{\t}{^{\mbox{\tiny\sf T}}}

\begin{document}

\begin{frontmatter}

\title{\LARGE Coupling effect and pole assignment in trajectory regulation of multi-agent systems}
\thanks[footnoteinfo]{
\newline
* Corresponding author.
}
\vspace{-.5cm}
\author[SWJU1]{Jilie Zhang}\ead{jilie0226@163.com},
\author[NC]{Zhiyong Chen}\ead{zhiyong.chen@newcastle.edu.au},
\author[HIT]{Hongwei Zhang}* \ead{hwzhang@hit.edu.cn},
\author[SWJU1]{Tao Feng} \ead{sunnyfengtao@163.com}

\address[SWJU1]{School of Information Science and Technology, Southwest Jiaotong University, Chengdu, Sichuan, 611756, P.R. China.}
\address[NC]{School of Electrical Engineering and Computing, The University of Newcastle, Callaghan, NSW 2308, Australia}
\address[HIT]{School of Mechanical Engineering and Automation, Harbin Institute of Technology, Shenzhen, Guangdong 518055, P.R. China.}
\vspace{-.5cm}

\begin{keyword}
Coupling effect; multi-agent system; pole assignment; trajectory regulation
\end{keyword}

\begin{abstract}
This paper revisits a well studied leader-following consensus problem of linear multi-agent systems, while aiming at follower nodes' transient performance.  Conventionally, when not all follower nodes have access to the leader's state information, distributed observers are designed to estimate the leader's state, and the observers are coupled via communication network. Then each follower node only needs to track its observer's state independently, without interacting with its neighbors. This paper deliberately introduces certain coupling effect among follower nodes, such that the follower nodes tend to converge to each other cooperatively on the way they converge to the leader. Moreover, by suitably designing the control law, the poles of follower nodes can be assigned as desired, and thus transient tracking performance can also be adjusted.
\end{abstract}

\end{frontmatter}


\vspace{-.2cm}
\section{Introduction}

\vspace{-.2cm}
Cooperative control of multi-agent systems (MASs) has been extensively investigated for the past two decades. It remains to gain increasing attention in the control community for its widespread applications in the areas of microgrids, unmanned aerial vehicles (UAVs), unmanned ground vehicles (UGVs), social networks, and so on. For a comprehensive literature review, readers are referred to some survey papers \cite{CaoYuRenChen2013Review,KnornChenMiddleton2016Overview,OhaParkAhn2015Survey} and references therein.

\vspace{-.1cm}
Among various problem formulations of multi-agent systems, a fundamental one is the cooperative tracking control problem, also known as leader-following consensus problem or trajectory regulation problem. In this scenario, there is one leader node and a group of follower nodes, and all the follower nodes are driven to track the trajectory of the leader node \cite{RenBook2008,ZhangLewisDas2011TAC-optimalsyn}. A salient feature of this problem is that only part of the follower nodes can or need to have access to the leader's information due to limited communication capability, and thus a distributed control algorithm is required.

\begin{figure}[h!]
	\begin{center}
		\includegraphics[width=8cm]{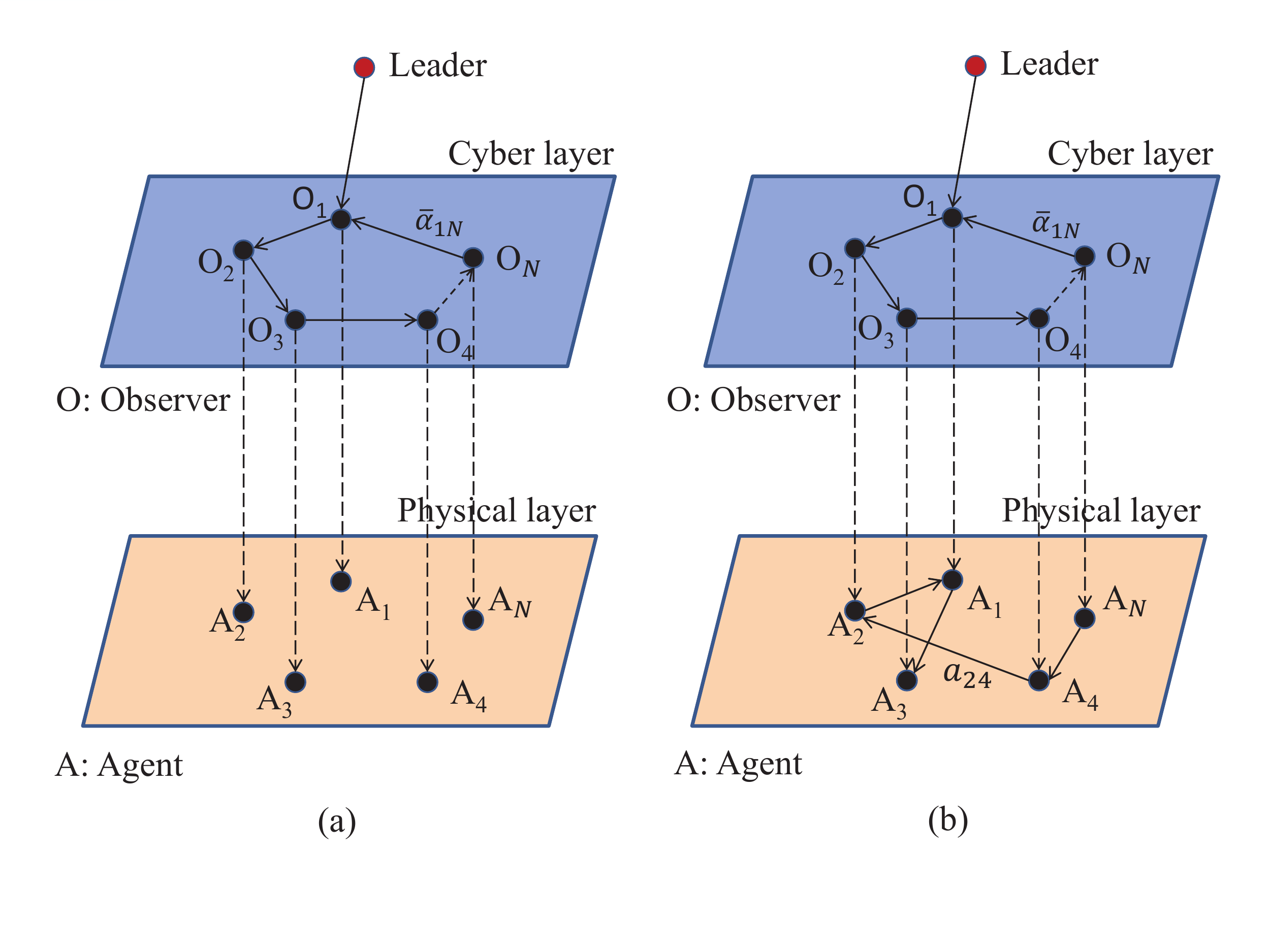}\vspace{-0.3in}
		\caption{Coupling in the cyber layer and physical layer}
		\label{fig:layer}
	\end{center}
\end{figure}

\vspace{-.15cm}
Such a problem has been well solved by introducing distributed observer design approach \cite{CaiH2017Auto-adaptive,Liang2019SMC-adaptive,SuYF2012TAC-output,ZhangLewisDas2011TAC-optimalsyn}, that is, each follower node maintains an observer, estimating the leader's state, and then each follower only need a local controller to drive its states to the leader's estimated state. In this framework, observer dynamics are coupled through communication network, while controllers are completely decoupled in a sense that each follower node only focuses on its own tracking task without considering its neighbors' information (see Fig. \ref{fig:layer}(a)). In other words, communication only exists in the cyber layer, instead of the physical layer.
 Although this control structure is very simple and intuitively understandable, it may result in a phenomenon that different follower nodes track the leader node independently, with some being uncoordinatedly faster/slower than others.  This is undesirable for some practical scenarios. For example,  in a competition task,
it is usually required that all agents, e.g., UAVs or UGVs, cooperatively arrive at their designated positions for a superior formation
simultaneously to avoid isolation and being vulnerable to enemy's attack.

\vspace{-.1cm}
It is also worth mentioning that almost all the existing works on distributed control of multi-agent systems focus on steady state collective behaviors, without considering transient performance of the group.  An exception is a line of research called prescribed performance control  \cite{Rovithakis2017TAC-decentr,YuMaZhang2019SMC-prescribed} or funnel control \cite{Shim2015CDC-funnel} of multi-agent systems. Both prescribed performance control and funnel control provide simple control laws that can restrict the profiles of the synchronization error within the prescribed error bounds. These control laws use neither the system dynamics information nor the graph topology information, and can be applied to linear systems for sure. However, their control laws consist of certain time functions, generating the prescribed performance bound profiles. This will lead to non-autonomous nonlinear closed-loop systems, whose performance depends on the initial time. Moreover, the initial values must be restricted within the prescribed bounds. Also, it is well known that  both transient and steady state performance of a linear system relies on the location of its poles. Then an interesting question is naturally raised that whether we can design a distributed control law such that the poles of all follower agents can be assigned as desired. This will provides an insight into the transient and steady state performance of multi-agent systems. To the best of our knowledge, this problem is still open.

\vspace{-.1cm}
Motivated by the above-mentioned statements, this paper aims to propose a novel distributed control law, which not only solves the trajectory regulation problem, but also considers the coupling effect between agents as well as desired pole assignment.
More specifically, we intend to deliberately introduce coupling effect for cooperative regulation, while avoiding over-coupling that may cause different issues, or even violation of system stability. Technically, compared with \cite{JL2020SMCA-State}, this paper gives a quantitative analysis of this trade-off in terms of stability,  dominant pole assignment, and fully system decoupling and hence pole assignment.
It is noted that a relevant work \cite{ygc2018} also introduces coupling effect between follower nodes,
but the design relies on the solution to an linear matrix inequality and each follower node must design its own coupling effect, and pole assignment is not considered.

\vspace{-.1cm}
The rest of this paper is organized as follows. The problem is formulated in Section \ref{s2}. Rigorous analysis of stability and pole assignment is provided in  Section \ref{Sect2}. A numerical example in Section \ref{s77} illustrates the effectiveness of the proposed algorithm. Section \ref{s7} concludes the paper. Technical lemmas are summarized in Appendix.

\vspace{-.1cm}
\textbf{Notations:}   The sets of real and complex numbers are denoted by  $\R$ and $\C$, respectively.
For a matrix $A$, $A\t$ represents it transpose and $A^*$ its conjugate transpose. For a symmetric real matrix $A$, $A>0$ ($A \geq 0$) means $A$ is positive definite (positive semi-definite).  The determinant of a matrix $A$ is denoted as det($A$). The real part of a complex number $m$ is denoted as $\rm Re [m]$. The set $\mathcal{N}$ is defined as $\{1,2,\dots,N\}$ and $\textbf{1}_N$ is the $N$-dimensional column vector whose elements are $1$.
The Kronecker product is denoted by the operator $\otimes$. For vectors $x_1, \cdots, x_n$,  $\mbox{col} (x_1, \cdots, x_n) =
[x\t_1, \cdots, x\t_n]\t$ represents the stacked vector.
\vspace{-.1cm}

\vspace{-.6cm}
\section{Problem Formulation}\label{s2}

\vspace{-.2cm}
The paper is concerned with control of a group of linear homogenous agents of the  dynamics described by
\vspace{-.2cm}
\begin{align}\label{sys1}
\dot{x}_i(t)=Ax_i(t)+Bu_i(t), \; i\in \mathcal{N}
\end{align}
where $x_i(t)\in \R^n$ is the state vector of the $i$-th agent and $u_i(t)\in \R^m$ the control input vector.
The control task is to regulate every state trajectory $x_i(t)$ to
the desired trajectory $x_0(t)\in \R^n$  described by the leader node
\vspace{-.2cm}
\begin{align}\label{sys10}
\dot{x}_0(t)=Ax_0(t),
\end{align}
\vspace{-.2cm}
that is,
\vspace{-.2cm}
\begin{align}\label{E3}
\lim_{t \to \infty}[x_i(t)-x_0(t)]=0, \; i\in \mathcal{N}.
\end{align}

\vspace{-.2cm}
The task becomes complicated when $x_0(t)$ is unaccessible for some agents. Nevertheless, such a task has
been well accomplished by using the framework of consensus of observers and trajectory regulation.
More specifically, an observer is established for each agent as follows
\begin{align}
\dot{\hat{x}}_{0i}= A\hat{x}_{0i} + \gamma_i (e), \; i\in \mathcal{N} \label{dhatx0i}
\end{align}
where $e=\mbox{col}(e_1,\cdots, e_N)$ with $e_i = \hat x_{0i} - x_0$ and the consensus protocol function $\gamma_i$ with $\gamma_i(0)=0$ is designed such that $\lim_{t \rightarrow \infty}e(t) =0$. Note that the design of $\gamma_i$ has been well studied in literature, and is not an issue considered in this paper. See the following remark for one choice of $\gamma_i$.
\begin{remark}
For homogenous MASs,  the design of consensus tracking protocol has been well studied \cite{SuYF2012TAC-output, ZhangLewisDas2011TAC-optimalsyn}, e.g., in a directed graph $\mathcal{G}$ represented by the network adjacent weight $\bar{\alpha}_{ij} \geq 0$, $i,j\in \mathcal{N}$,
\vspace{-.2cm}
\begin{align}
\gamma_i (e)&=   \mu \left(\sum_{j=1}^N \bar{\alpha}_{ij}(e_{ j}(t)-e_{ i}(t)) - g_{i}e_i(t) \right) \nonumber\\
&=  {\mu}\left(\sum_{j=1}^N  \bar{\alpha}_{ij}(\hat{x}_{0 j}(t)-\hat{x}_{0 i}(t))+g_{i}({x}_0(t)-\hat{x}_{0i}(t))\right)\label{gamma}
\end{align}
with  $\mu$ sufficiently large \cite{SuYF2012TAC-output} (or $\mu$ can be designed in a decoupled manner $\mu=cK$ as in \cite{ZhangLewisDas2011TAC-optimalsyn}) and $g_i\neq 0$ for $i \in \mathcal{S_R}$,  $g_i=0$ for $i \notin \mathcal{S_R}$,
where $\mathcal{S_R} \subset \mathcal{N}$ is the set of  nodes that have  direct access to the reference signal $x_0(t)$.
In this setting, if every follower node $i\in \mathcal{N}$ is reachable in $\mathcal{G}$ from the leader node, then one has $\lim_{t \to \infty}e(t) =0$.
\end{remark}

\vspace{-.3cm}
The time variable $(t)$ is ignored in \eqref{dhatx0i} and the sequel for conciseness.
Next, let $\eta=\mbox{col}(\eta_1,\cdots, \eta_N)$ with
\vspace{-.2cm}
$$\eta_i =x_i -  \hat{x}_{0i}.$$
A trajectory regulator
\vspace{-.1cm}
\begin{align}
u_i=- d_i F  \eta_i,~~ i\in \mathcal{N} \label{contr0}
\end{align}
with a feedback matrix $F\in \R^{m\times n}$ and a scalar gain $d_i >0$
is designed such that $\lim_{t \to \infty}\eta(t) =0$, that is,
the state $x_{i}(t)$ is regulated to the observer's state $\hat{x}_{0i}(t)$.
As the closed-loop system is
\vspace{-.1cm}
\begin{align*}
\dot{\eta}_i =(A -d_i B F )\eta_i  - \gamma_i (e), \; i\in \mathcal{N},
\end{align*}
it suffices to pick $F$ and $d_i$ such that $A -d_i B F$ is Hurwitz.
From above, the task \eqref{E3} is accomplished by combining $\lim_{t \rightarrow \infty}e(t) =0$
and $\lim_{t \to \infty}\eta(t) =0$.

Now, the main focus of this paper is on the design of regulator for $\lim_{t \rightarrow \infty}\eta(t) =0$.
It is obvious that the simple regulator \eqref{contr0} is designed separately for each individual agent. Its simplicity has
independent interest in many scenarios. However, in some other scenarios, researchers are also interested
in more sophisticated interactions among agents when they obey the trajectory regulation protocol. For instance,
individual regulation does not work when the agents intend to converge to each other before to the desired trajectory.
For this purpose, we {\it deliberately}  introduce coupling effect $a_{ij}$ among agents in the physical layer (see Fig. \ref{fig:layer}(b)), which results in
\vspace{-.2cm}
\begin{align}
u_i=&F \left( \sum_{j=1}^N a_{ij} \eta_j -d_i \eta_i \right),\;i\in \mathcal{N}, \label{contr-c}
\end{align}
where $a_{ij}\geq 0$, $i,j\in \mathcal{N}$ represents the coupling weight between agent $i$ and $j$ in the physical layer, with
$a_{ii}=0$.  Note that $\bar{\alpha}_{ij}$ denotes the coupling among observers in the cyber layer (see Fig. \ref{fig:layer}).
It is worth mentioning that the leader model \eqref{sys10} has the role in the design
of the observer network in \eqref{gamma}, but it has no direct influence on the
deliberate coupling network in \eqref{contr-c}. In other words, the main results given in this paper,
based on the analysis of \eqref{contr-c}, may apply not only to a leader-following topology, but also to
more general leaderless topologies.

The closed-loop MASs composed of \eqref{sys1},  \eqref{dhatx0i}, and \eqref{contr-c}
can be written as
\vspace{-.2cm}
\begin{align*}
\dot{\eta}_i =A\eta_i+BF \left( \sum_{j=1}^N a_{ij} \eta_j -d_i \eta_i \right)  - \gamma_i (e),\;i\in \mathcal{N}.
\end{align*}
Let ${\mathcal A}$
be a matrix whose $(i,j)$-entry is $a_{ij}$,  i.e., ${\mathcal A} = [a_{ij}]$,  ${\mathcal D} =\mbox{diag} (d_1,\cdots, d_N)$ be a diagonal matrix,
and ${\mathcal M}= {\mathcal D}-{\mathcal A}$.
The closed-loop system can rewritten in a compact form
\begin{align*}
\dot{\eta} &= [(I_N\otimes A)-({\mathcal M}\otimes BF)] \eta  - \gamma (e) \nonumber\\
&= (\bar{A}-\bar{B}\bar{\mathcal M}\bar{F})\eta  - \gamma (e)
\end{align*}
where $\gamma (e) =\mbox{col}(\gamma_1(e),\cdots, \gamma_N(e))$ and
\begin{align*}
\bar{A}=I_N\otimes A,\; \bar{B}=I_N\otimes B,\; \bar{F}=I_N\otimes F, \;
\bar{\mathcal M} ={\mathcal M}\otimes I_m.
\end{align*}
Alternatively,  let $\xi=\mbox{col}(\xi_1,\cdots, \xi_N)$ with
$\xi_i =x_i -  {x}_{0}$. Note that $e=\xi-\eta$.
A straightforward computation shows
\vspace{-.2cm}
\begin{align*}
\dot{\xi} =( \bar{A} -\bar{B}\bar{\mathcal M}\bar{F}) \xi +\bar{B}\bar{\mathcal M}\bar{F} e.
\end{align*}

\vspace{-.2cm}
With $\lim_{t \to \infty}e(t) =0$, it suffices to establish a stable system
\begin{align}\label{deta}
\dot{\eta} = (\bar{A}-\bar{B}\bar{\mathcal M}\bar{F})\eta,
\end{align}
or equivalently,
\begin{align}\label{dxi}
\dot{\xi} =( \bar{A} -\bar{B}\bar{\mathcal M}\bar{F}) \xi
\end{align}
for accomplishment of the task \eqref{E3}.

With the above development, this paper aims to analyze the coupling effect of $\mathcal A$,
in comparison with $\mathcal D$. First of all, it should be noted that
(i) The diagonal entries of ${\mathcal M}$ are represented by $\mathcal D$ and the off-diagonal entries
by $-\mathcal A$.  If $d_i = \sum_{j=1}^Na_{ij}$, $i\in\mathcal N$,  the matrix ${\mathcal M}$ becomes a Laplacian.
But in general, the matrix ${\mathcal M}$ is adjusted as a diagonally dominant matrix for our purpose
with a more significant $\mathcal D$.
(ii) The graph for the adjacency matrix $\mathcal A$ is not necessarily connected for the purpose
of consensus. Especially, the consensus can be achieved by the controller with ${\mathcal A}=0$, i.e.,
\eqref{contr0}.

In this paper, on one hand, we intend to introduce more significant coupling effect $a_{ij}$ on the physical layer for cooperative regulation.
On the other hand, over-coupling may cause different issues even violation of system stability. So,
the main objective of this paper is to conduct
more specific analysis of the coupling effect of $\mathcal A$ relative to $\mathcal D$ from the following three aspects,
which apply different upper boundary conditions on $a_{ij}$. These conditions can be regarded as design criteria in selecting the strength of $\mathcal A$ for \eqref{deta} or \eqref{dxi}.

  {\it
 \begin{itemize}

 \item[(1)] The condition on $\mathcal A$  relative to $\mathcal D$ for stability
 of \eqref{deta} or \eqref{dxi}.

 \item[(2)] A less significant $\mathcal A$ and a more significant $\mathcal D$ allow dominant
pole assignment of \eqref{deta} or \eqref{dxi}.

 \item[(3)] The condition on $\mathcal A$  relative to $\mathcal D$ for a
 fully decoupled system matrix and hence pole assignment.
  \end{itemize}
 }

\section{Main Results} \label{Sect2}
\vspace{-.1cm}
The three aspects listed in the aforementioned main objective are addressed in this section in order.
In the subsequent presentation, we will present the results on \eqref{dxi}.

\subsection{Stability Condition}

Each agent uses its own feedback through $d_i$ and others via $a_{ij}$.
Intuitively, it needs more significant feedback from itself for stability.
To explicitly characterize the relative significance, we introduce the condition under which
${\mathcal M} + {\mathcal M}\t$ is positive definite. This property plays an important role in proving  the stability of  \eqref{dxi}.

\vspace{-.1cm}
\begin{lemma} \label{L2}
For ${\mathcal M}= {\mathcal D}-{\mathcal A}$,  if
\begin{align}\label{di-stab}
d_i>\frac{\sum_{j=1}^N (a_{ij} + a_{ji})}{2 },\; i \in \mathcal{N},
\end{align}
then
${\mathcal M} + {\mathcal M}\t>0$.
\end{lemma}
\vspace{-.3cm}
\bproof
According to the construction of ${\mathcal M}= {\mathcal D}-{\mathcal A}$, the diagonal entry of the $i$-th row of
${\mathcal M} + {\mathcal M}\t$ is $o_i= 2 d_i$ and the sum of the absolute values of the off-diagonal entries in the $i$-th row is
$r_i= {\sum_{j=1}^N (a_{ij} + a_{ji}})$. As a result, one can define a Ger\v{s}gorin disc, denoted by $D(o_i, r_i)$, which is centered at $o_i$
with radius $r_i$.

Under the condition \eqref{di-stab}, i.e., $o_i > r_i$, $ i \in \mathcal{N}$,
all the  Ger\v{s}gorin discs are located in the open right half plane in the complex plane. Therefore, by Ger\v{s}gorin theorem, all eigenvalues of $ {\mathcal M} + {\mathcal M}\t$ are positive real numbers, which concludes
${\mathcal M} + {\mathcal M}\t>0$.
\eproof
\vspace{-.3cm}
\begin{figure}[H]
 \centering
  \includegraphics[width=7cm]{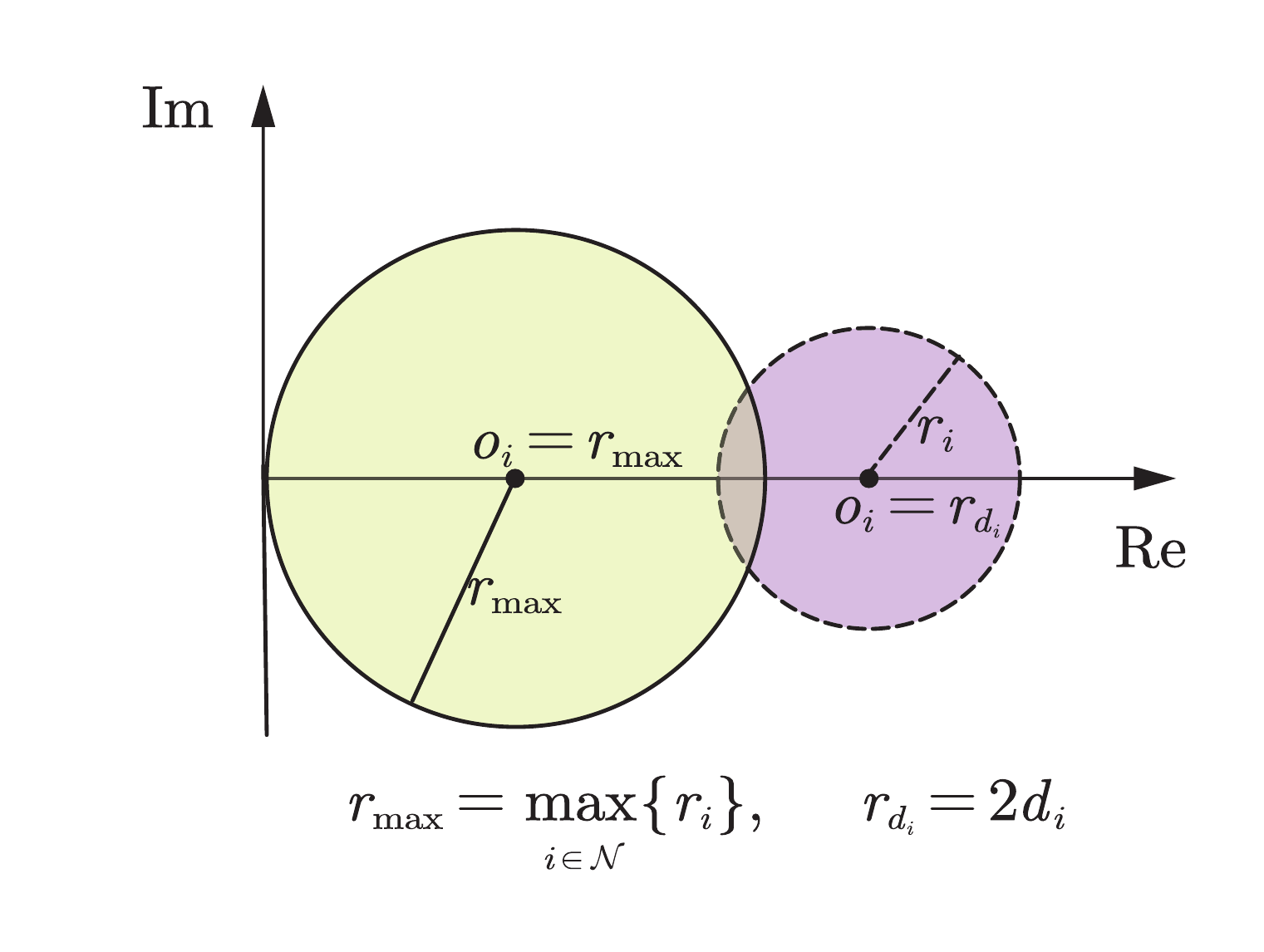}\vspace{-.3cm}
\caption{Illustration of Ger\v{s}gorin discs of ${\mathcal M} + {\mathcal M}\t$. }\label{f1}
\end{figure}

\begin{remark}
Lemma \ref{L2} can be intuitively shown by Fig.\ref{f1}. Let $r_{\max} =\max_{i\in \mathcal N}\{r_i\}$. If
$d_i=\sum_{j=1}^N(a_{ij} + a_{ji}) /2 ,\; i \in \mathcal{N}$, ${\mathcal M} + {\mathcal M}\t$ becomes a symmetric Laplacian.
All Ger\v{s}gorin discs are located inside of the largest one $D(r_{\max}, r_{\max})$, illustrated by the solid circle,
and go through the origin. However, under the condition \eqref{di-stab}, all the Ger\v{s}gorin discs are strictly shifted to the right, not going through the origin any more,  illustrated by the dotted circle.
\end{remark}

\begin{figure}[h!]
  \centering
  \includegraphics[width=2in]{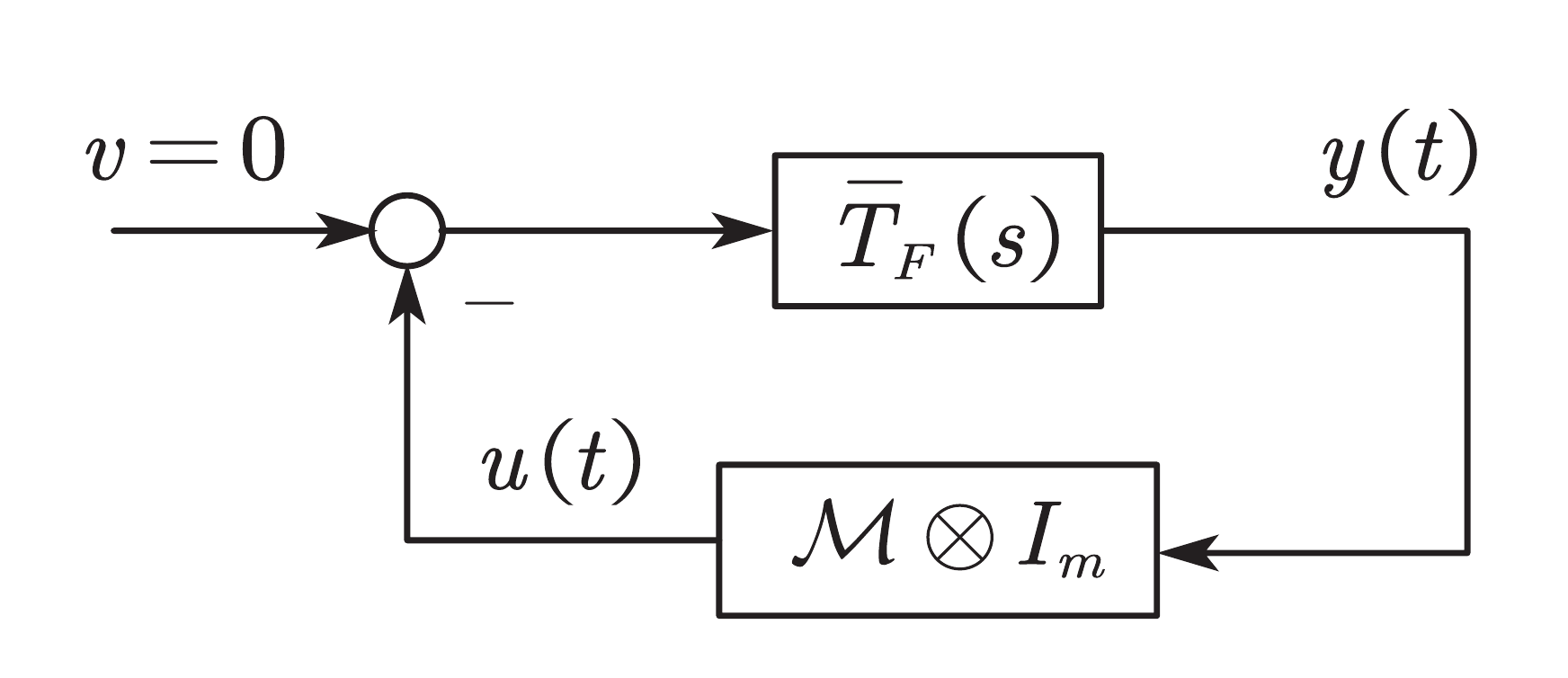}
  \caption{Schematic diagram of an equivalent state feedback system. }\label{f2}
\end{figure}

\vspace{-.1cm}
Note that the stability of the system  \eqref{dxi} can equivalently represented by
a feedback structure illustrated  in Fig.~\ref{f2} in frequency domain. In other words,
\eqref{dxi} can be rewritten as follows, with $v=0$,
\vspace{-.1cm}
\begin{align}\label{E23}
\dot{\xi} =&\bar{A}{\xi}+\bar{B}u,  \nonumber\\
y=&\bar{F}\xi, \nonumber\\
u=&-({\mathcal M}\otimes I_m){y} + v.
\end{align}
Here $v$ is an auxiliary input introduced for the convenience of presentation later.
Now, the main result is stated below.

\vspace{-.1cm}
\vspace{-.1cm}
\begin{assumption}\label{A1}
The matrix pair $(A,B)$ is stabilizable.
\end{assumption}

\vspace{-.2cm}
\begin{theorem}\label{T1}
Under Assumption \ref{A1} and the condition \eqref{di-stab}, there exists $F$ such that the system \eqref{dxi} is stable.
\end{theorem}

\vspace{-.2cm}
\bproof A constructive proof is given as follows. First of all,
Lemma \ref{L2} shows that ${\mathcal M} + {\mathcal M}\t>0$ under the condition \eqref{di-stab}.
There exists a scalar $k>0$ such that
$ {\mathcal M} + {\mathcal M}\t > 2 k I$. Then,
$Z+Z\t>0$ for $Z = ({\mathcal M} - k I )/k$. That is,
${\mathcal M} =k (Z +I)$.
 Now, the system \eqref{E23} can be rewritten as
 \vspace{-.1cm}
\begin{align}
\dot{\xi} =&(\bar{A} -k \bar{B} \bar{F} ) \xi  +\bar{B} \bar u, \nonumber\\
  y=&\bar{F}\xi, \nonumber\\
\bar u=& - k (Z \otimes I_m){y} +v. \nonumber
\end{align}

\vspace{-.3cm}
Under Assumption \ref{A1}, there exists a solution $X_{C}> 0$ to the control algebraic Riccati equation
\begin{align*}
A\t  X_{C} + X_{C}A- k^2 X_{C}BB\t X_{C} + Q_{C} = 0
\end{align*}
with $Q_{C}>0$.
Let $F = - k B\t X_{C}$.
By Lemma~\ref{F2}, the transfer function matrix $T_{F}(s) = F(sI-A + kBF)^{-1} kB$
is positive real, so is  $\overline{T}_F(s) =I_N\otimes T_{F}(s)$.
In other words,
\begin{align}
T_F(s)+T_{F}^*(s) \geq 0, \; \forall\ {\rm Re}[s]>0.
\label{TFs}
\end{align}

\vspace{-.2cm}
It is easy to see that the transfer function from $v$ to $y$ is
 \begin{align} \label{transfer}
I_N \otimes \frac{T_F(s) }{I+ k T_F(s) Z}
\end{align}
whose poles are determined by  $\det(I + k {T}_F(s) Z) =0$.
One has
\begin{align*}
\det[I + k {T}_F(s) Z]  \neq 0, \; \forall\ {\rm Re}[s] \geq 0, \end{align*}
 due to \eqref{TFs} and $Z+Z\t>0$, using Lemma~\ref{FF}. It means that all the poles
 of \eqref{transfer} have negative real parts.
The stability of \eqref{dxi}  is thus proved.
\eproof

\begin{remark} A uniform selection of $d_i$, i.e.,
\begin{align*}
d_i = d > \max_{i\in \mathcal N} \left\{ \frac{\sum_{j=1}^N (a_{ij} + a_{ji})}{2 } \right\},
\end{align*}
can be used to satisfy the condition  \eqref{di-stab}. It simplifies the protocol \eqref{contr-c} with less parameters.
Condition  \eqref{di-stab} requires sufficiently large $d_i$.
However, high gain of a controller may lead to some practical issues.
To prevent the high gain issue, in practice, we can choose $d_i$ to the desired level and then scale $a_{ij}$ down to match the condition  \eqref{di-stab}.
\end{remark}

\subsection{Pole Assignment}\label{s66}

In this subsection, we aim to show that the system poles can be assigned when
$\mathcal D$ is sufficiently large and $\mathcal A$ sufficiently small.
In particular, the transient characteristics of the regulation behavior can be adjusted by pole assignment.
We will use the inverse optimal linear quadratic regulator (LQR) technique to assign poles.
Here we assume $\mbox{rank}(B)=m$ and the matrices $A, B$ take the following form
\begin{align}
A=\left[\begin{array}{cc}
A_{11} & A_{12} \\
A_{21} & A_{22}
\end{array}
\right],\;
B=\left[\begin{array}{c}
0\\
B_2
\end{array}
\right],\; (\det (B_2) \neq 0) \label{AB2}
\end{align}
where  $A_{22}$ and $B_2$ have the same dimension.
Define the LQR problem
\begin{align}
\dot z &= A z + B w \nonumber \\
J &= \int_{0}^\infty (z\t (t) Q z(t) + w\t(t) R w(t)) dt \label{LQ}
\end{align}
with $Q=Q\t \geq 0$ and $R=R\t >0$ for the convenience of presentation.
Now, the main result is stated in the following theorem.

\begin{theorem}\label{T5}
Consider the system \eqref{dxi} of the structure \eqref{AB2} under Assumption \ref{A1}.
Let
\begin{align*}
 F = B_2^{-1} V^{-1} \Sigma V[K\;\; I] ,\; \Sigma = \mbox{diag} (\sigma_1, \cdots, \sigma_m),
\end{align*}
where $\sigma_1, \cdots, \sigma_m>0$ and
$V, \Sigma, K$ are selected such that $w = -F z$ is a solution to the LQR problem \eqref{LQ}.
In particular, the eigenvalues of $A_{11} -A_{12} K$ has the specified stable  eigenvalues
$\{s_{1}, \cdots, s_{n-m}\}$, none of which is an eigenvalue of $A_{11}$.
If  \begin{align} \label{diA}
d_i \to \infty, \; i\in \mathcal{N},\;  \|{\mathcal A}\| \to 0,
\end{align}
the system \eqref{dxi}  has the following eigenvalue distribution:
\begin{itemize}

\item[(i)] there are $(n-m)N$ eigenvalues $\lambda_i^j$ of the form $\lambda_i^j \to s_j$,
$j=1,\cdots, n-m$, $i\in \mathcal N$; and

\item[(ii)] there are $mN$ eigenvalues $\lambda_i^{j+n-m}$ of the form $\lambda_i^{j+n-m} \to -d_i \sigma_j$,
$j=1,\cdots, m$, $i\in \mathcal N$.
 \end{itemize}

 \end{theorem}
\vspace{-.2cm}
\bproof
First of all, it is noted that the system \eqref{dxi} can be rewritten as
\begin{align*}
\dot{\xi} =( \bar{A} -\bar{B}\bar{\mathcal D}\bar{F}) \xi  -\bar{B}\bar{\mathcal A}\bar{F} \xi.
\end{align*}
where $\bar{\mathcal D} ={\mathcal D}\otimes I_m$ and $\bar{\mathcal A} ={\mathcal A}\otimes I_m$.
For the eigenvalues of a matrix continuously depend on its parameter variation, the
eigenvalues $\lambda_i^j$, $j=1,\cdots, n$, $i\in \mathcal N$, approach those of $\bar{A} -\bar{B}\bar{\mathcal D}\bar{F}$
as $\|{\mathcal A}\| \to 0$. Furthermore, the eigenvalues of $\bar{A} -\bar{B}\bar{\mathcal D}\bar{F}$
are those of $A -d_i B   F$.

For every $A -d_i B   F$,  $i \in \mathcal N$, as $d_i \to \infty$,
there are $n-m$ eigenvalues $\lambda_i^j$ of the form $\lambda_i^j \to s_j$,
$j=1,\cdots, n-m$, and $m$ eigenvalues $\lambda_i^{j+n-m}$ of the form $\lambda_i^{j+n-m} \to -d_i \sigma_j$,
$j=1,\cdots, m$, by Lemma~\ref{pole}. The completes the proof. \eproof

\begin{remark}
The theorem shows that the eigenvalues of the closed-loop system approach arbitrarily specified
stable poles $s_j$ for $j=1,\dots,n-m$ and other stable poles $-d_i\sigma_j$ when all $d_i$ are sufficiently large and all $a_{ij}$ are sufficiently small. In particular, as all $d_i$ are sufficiently large, $ \rm Re [s_j] \gg -d_i\sigma_{j}$. In other words, all specified
$s_j$ are the dominant poles, enforcing the plant the specified transient characteristics.
 \end{remark}

\vspace{-.3cm}
\subsection{Fully Decoupling Condition}\label{s5}
\vspace{-.1cm}
In the previous subsection, it is proved that the dominant poles of
 the closed-loop system can be placed as desired
such that the transient characteristics can be satisfied, provided that $\|{\mathcal A}\| \to 0$.
It obviously contradicts to the main motivation of adding coupling effect to regulation,
if we simply remove the coupling by letting ${\mathcal A} =0$. In this subsection,
we will further study an explicit condition on the size of $a_{ij}$ under which
the system can be fully decoupled and the pole assignment technique can be
applied without the assumption of $\|{\mathcal A}\| \to 0$.

We first give a lemma for diagonalization of the coupling matrix ${\mathcal M}$, followed by the main theorem.

\begin{figure}[h!]
  \centering
  \includegraphics[width=8cm]{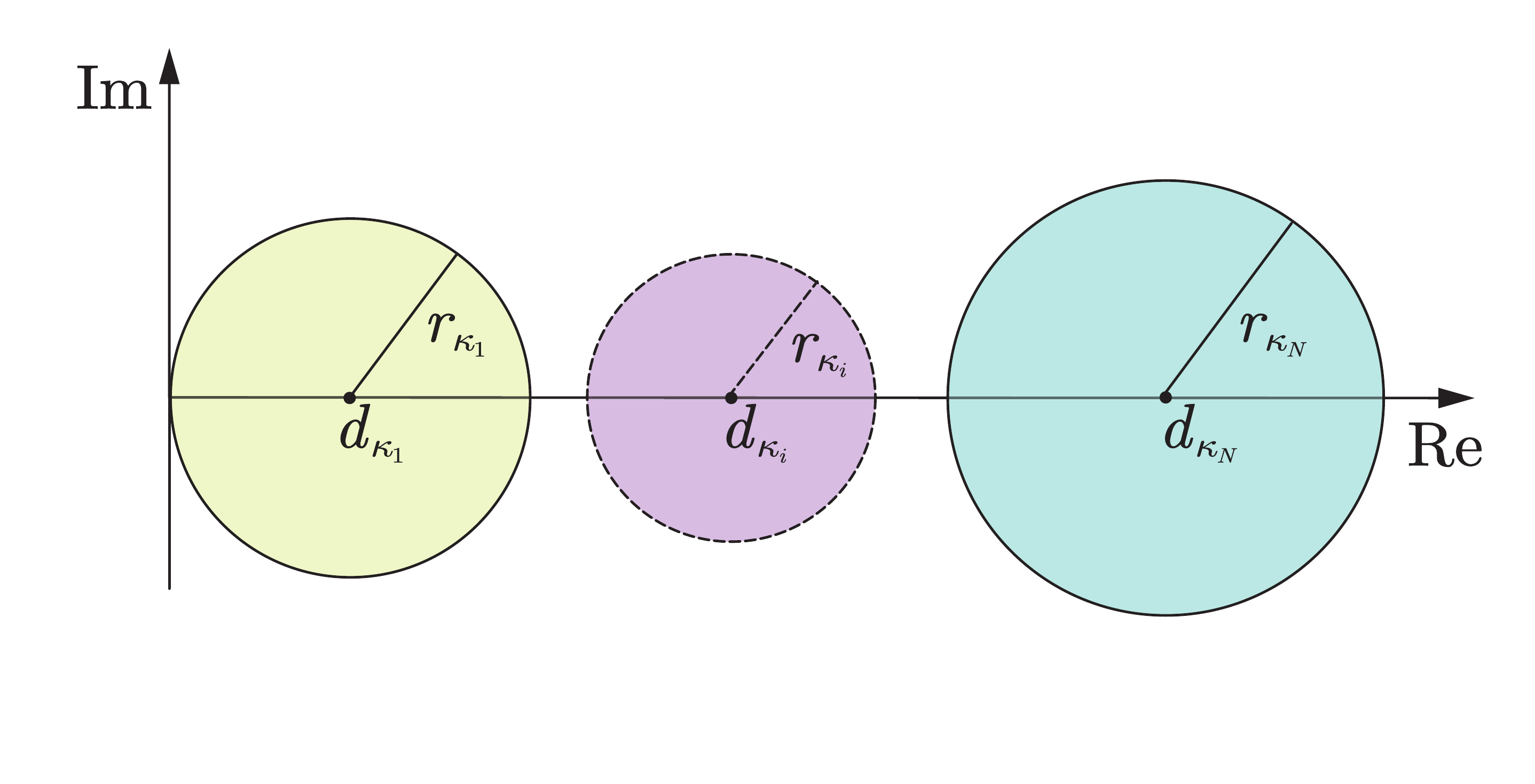}\vspace{-.3cm}
  \caption{Illustration of non-intersected Ger\v{s}gorin discs.}\label{f22}
\end{figure}

\begin{lemma}\label{l4}
For ${\mathcal M}={\mathcal D}-{\mathcal A}$,
if, for a sequence $\{\kappa_n\}$ with $\{ \kappa_1, \cdots, \kappa_N\} =\mathcal N$ and
\begin{align}
d_{\kappa_1} &> r_{\kappa_1},\; \notag \\
d_{\kappa_{i+1}}-d_{\kappa_{i}} &> r_{\kappa_{i+1}}+r_{\kappa_{i}}, \; i=1,\ldots, N-1, \label{ddiag}
\end{align}
where $r_i= \sum_{j=1}^N   a_{ij}$,  then ${\mathcal M}$ is diagonalizable.
\end{lemma}
\vspace{-.5cm}
\bproof The diagonal entry of the $i$-th row of
${\mathcal M}$ is $d_i$ and the sum of the absolute values of the off-diagonal entries in the $i$-th row is
$r_i= \sum_{j=1}^N   a_{ij}$. As a result, one can define a Ger\v{s}gorin disc centred at $d_i$
with radius $r_i$, denoted $D(d_i, r_i)$.
Under the condition \eqref{ddiag}, the Ger\v{s}gorin discs are
$D(d_{\kappa_i}, r_{\kappa_i})$ in the order of $i=1,\cdots, N$,
from left to right, and they do not intersect (see Fig. \ref{f22}).  Therefore, by Ger\v{s}gorin theorem, all eigenvalues of $ {\mathcal M}$ are distinct positive real numbers. Thus ${\mathcal M}$ is diagonalizable.
\eproof
\vspace{-.1cm}
\vspace{-.2cm}
 \begin{theorem}\label{T6}
Consider the system \eqref{dxi} of the structure \eqref{AB2} under Assumption \ref{A1}
and the condition \eqref{ddiag}.  Theorem~\ref{T5} holds with \eqref{diA} replaced by
\vspace{-.1cm}
\begin{align*}
d_i \to \infty, \; i\in \mathcal{N}.
\end{align*}
 \end{theorem}
 \vspace{-.4cm}
\bproof By Lemma~\ref{l4}, there exists a nonsingular matrix $T$ such that
\begin{align*}
 \vspace{-.1cm}
{\mathcal P} = T {\mathcal M} T^{-1},\; {\mathcal P} = \mbox{diag} (p_1,\cdots, p_N)
\end{align*}
and $p_i \to \infty, \; i\in \mathcal{N}$ as $d_i \to \infty, \; i\in \mathcal{N}$.

 \vspace{-.1cm}
With $\zeta = (T \otimes I_m) \xi$,  the system \eqref{dxi} is equivalent to
 \vspace{-.1cm}
\begin{align*}
\dot{\zeta} = ( \bar{A} -  \bar{B}({\mathcal P} \otimes I_m )\bar{F} ) \zeta
\end{align*}
whose poles are determined by the eigenvalues of $A -p_iBF$.
The remaining proof follows that of Theorem~\ref{T5}.
\eproof

 \vspace{-.1cm}
\section{Numerical Example}\label{s77}
 \vspace{-.1cm}
Consider a group of $N=5$ linear homogenous agents of the  dynamics described by \eqref{sys1} with
 \begin{align*}
 A=\left[
\begin{array}{cc}
 0 & -0.5 \\
0 & 0 \\
 \end{array}
 \right], \;
 B=\left[\begin{array}{c}
 0  \\
 1 \\
 \end{array}
 \right].
 \end{align*}
In the simulation, let $F=\left[
            -0.3660,\; 0.9306
        \right]$ as in  Theorem~\ref{T5},
and pick $\mathcal D =\mbox{diag}\{2,1.95,2.5,2.9,3.4\} \rho$ and $\mathcal A= [\varepsilon a_{ij}]$ with $a_{12}=a_{23}=a_{34}=a_{45}=a_{51}=1$ and $a_{ij}=0$ otherwise,
for two parameters $\rho$, $\varepsilon >0$.
The results are demonstrated in terms of
  $x_i = \xi_i +  {x}_{0}$ with $x_0(t)$ governed by \eqref{sys10} and
 $\xi(t)$ by \eqref{dxi}, where the specific behavior of $e(t)$ approaching zero is ignored.

 Denote $\xi_i = [\xi_{i1}, \xi_{i2}]\t$ and $x_i = [x_{i1}, x_{i2}]\t$. Also, denote
 $\theta = [x_{11},\cdots, x_{51}]\t$ as the vector consisting of the first elements of the five agent states.
 The signal $\phi = \max\{ |\xi_{11}|,\cdots, |\xi_{51}|\}$ represents the regulation error
 and  $\psi= \max\{ \xi_{11},\cdots, \xi_{51}\} -\min\{ \xi_{11},\cdots, \xi_{51}\}$
 the difference among the agents.
\vspace{-.3cm}
 \begin{figure}[H]
  \centering
\hspace{-10mm}
  \includegraphics[width=8cm]{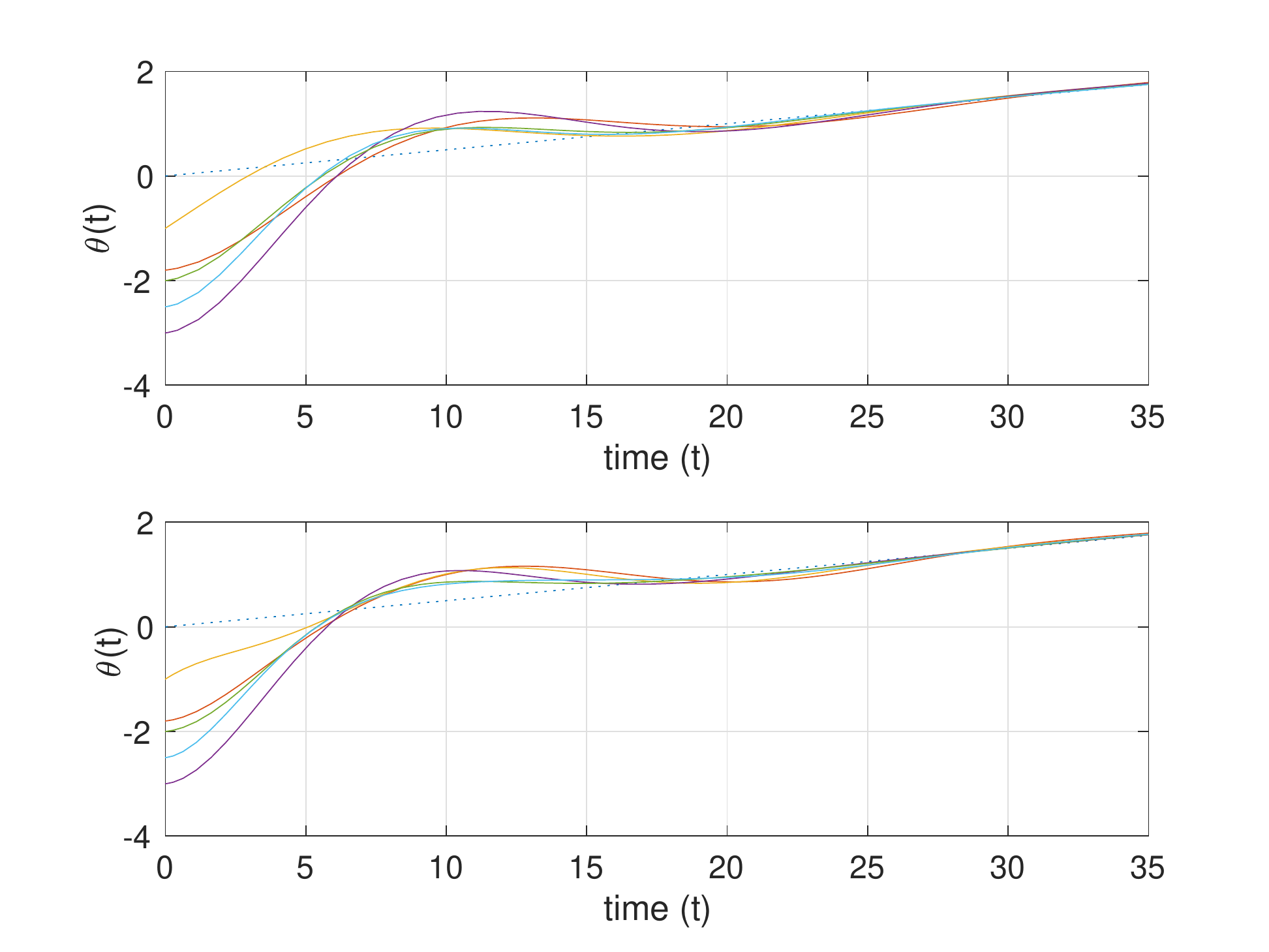} \vspace{-.2cm}
  \caption{Profile of five agents achieving synchronization; top: $\rho =0.12$, $\varepsilon =0$;
  bottom:  $\rho =0.2$, $\varepsilon =0.2$.  }\label{fig1}
\hspace{-10mm}
  \includegraphics[width=8cm]{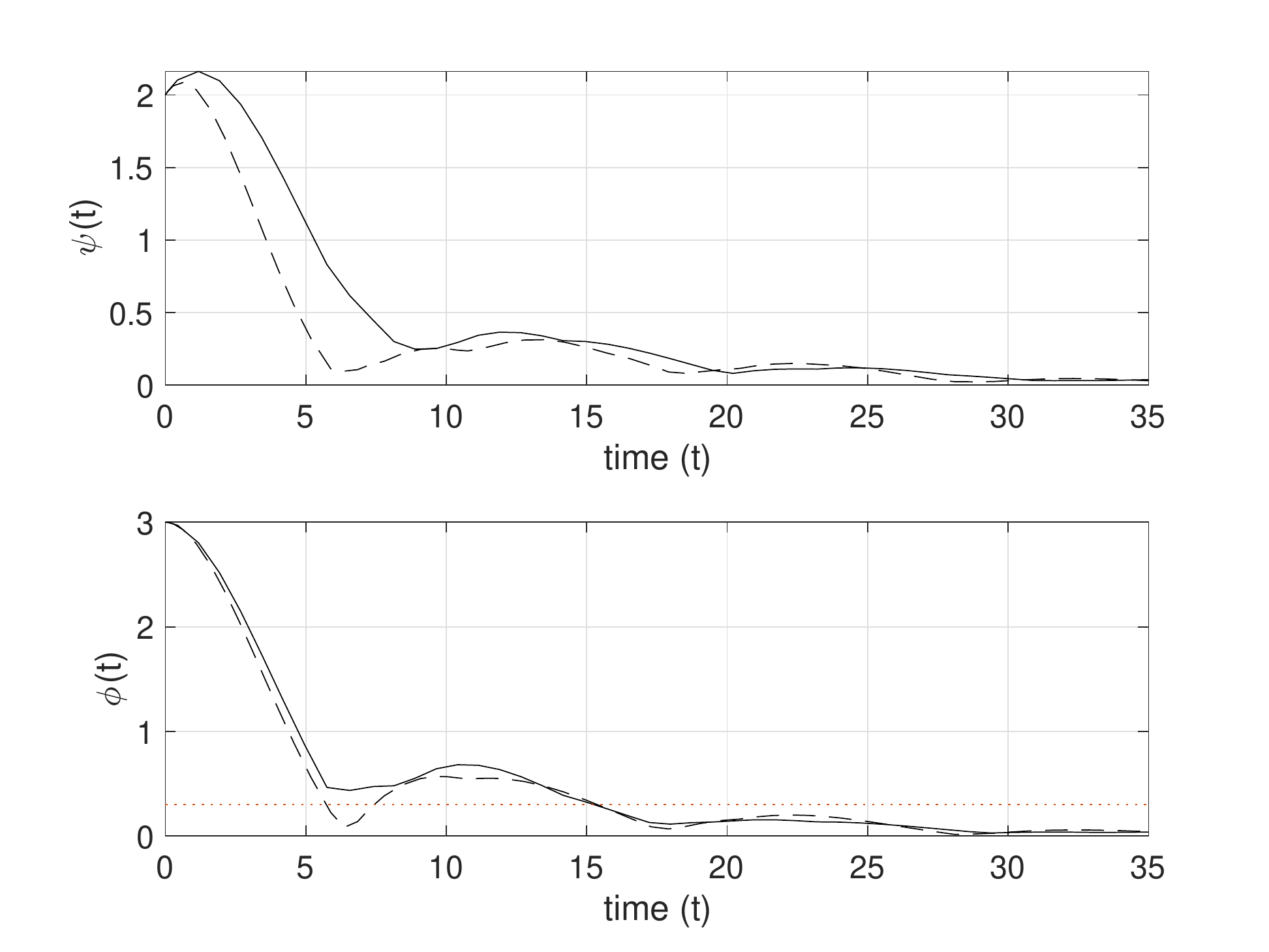} \vspace{-.2cm}
  \caption{Profile of agent difference $\psi(t)$ (top) and regulation error $\phi(t)$ (bottom); solid lines: $\rho =0.12$, $\varepsilon =0$;
dashed lines:  $\rho =0.2$, $\varepsilon =0.2$.}\label{fig2}
\end{figure}
\vspace{-.4cm}
In the first case, we select $\rho =0.12$ and $\varepsilon =0$.
As shown in Fig.~\ref{fig1},  the agents reach consensus while their transient processes are independent and do not influence each other. In other words, the agents do not demonstrate a cooperative behavior.
The profiles of the agent difference and the regulation error are shown in Fig.~\ref{fig2}.

\vspace{-.1cm}
In the second case, we select $\rho =0.2$ and $\varepsilon =0.2$ as comparison.
All the agent states also converge a consensus on the same reference signal.
With $\varepsilon =0.2$, the agents have cooperation  before achieving consensus.
 In the simulation, we pick the  parameters such that the regulation error achieves zero with
 the same transient performance. As shown in the bottom graph of Fig.~\ref{fig2}, the error reduces from 3 to 0.3 (by $90\%$)
 at $t=15$ for both two cases. It is obvious that,
 the difference among agents in  the case with $\varepsilon =0.2$
 is significantly less than that in the case with $\varepsilon =0$.
 More specifically, it takes $t=4.5$ for the difference among agents to reduce below $0.5$ in the former case
 while it takes $t=7.0$ in the latter case, as shown in the top graph of Fig.~\ref{fig2}.

\begin{figure}[H]
    \centering
\hspace{-10mm}
  \includegraphics[width=8cm]{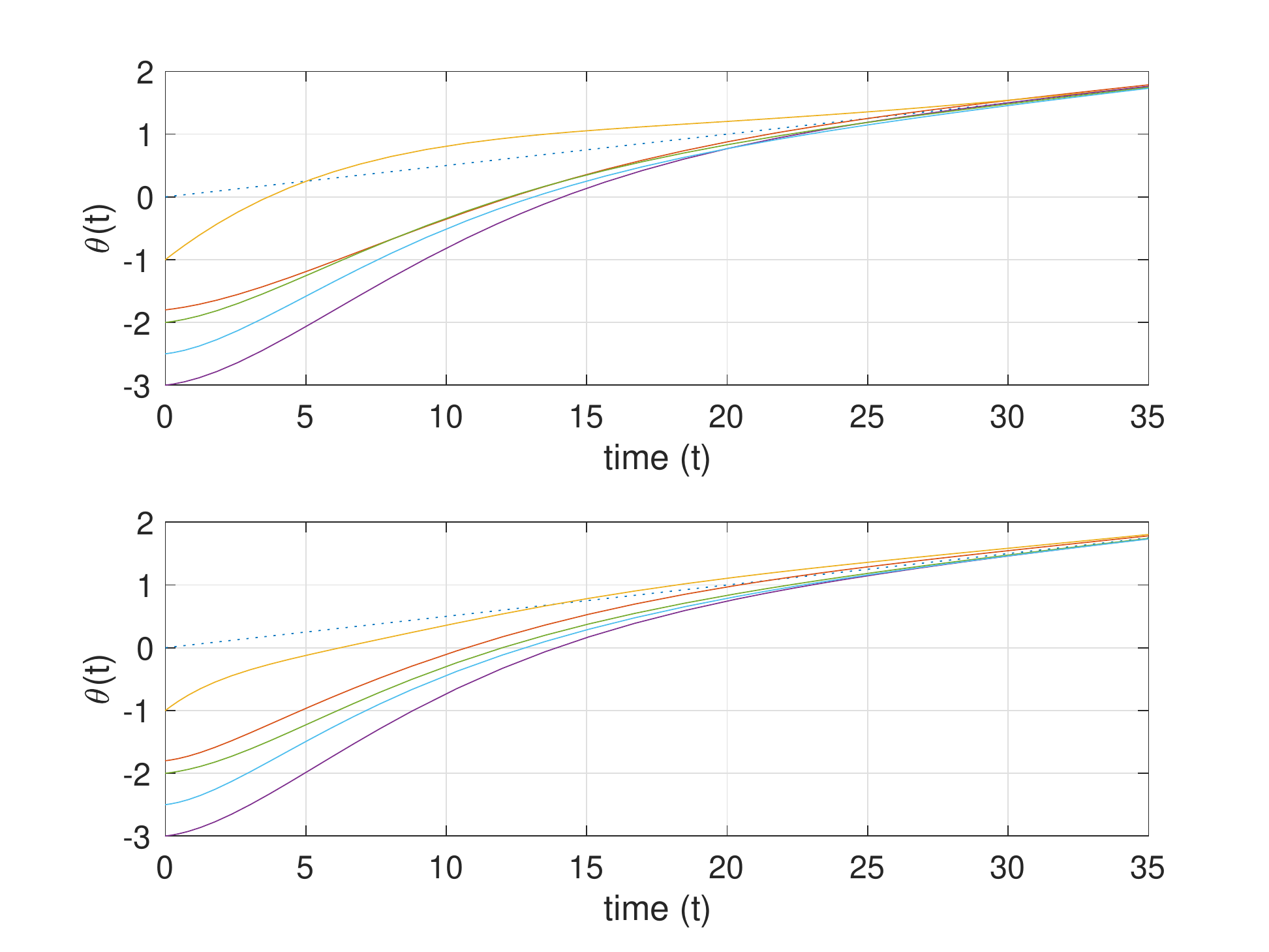}\vspace{-.3cm}
  \caption{Profile of five agents achieving synchronization; top: $\rho =0.12$, $\varepsilon =0$;
  bottom:  $\rho =0.2$, $\varepsilon =0.2$.  }\label{fig3}
\hspace{-10mm}
  \includegraphics[width=8cm]{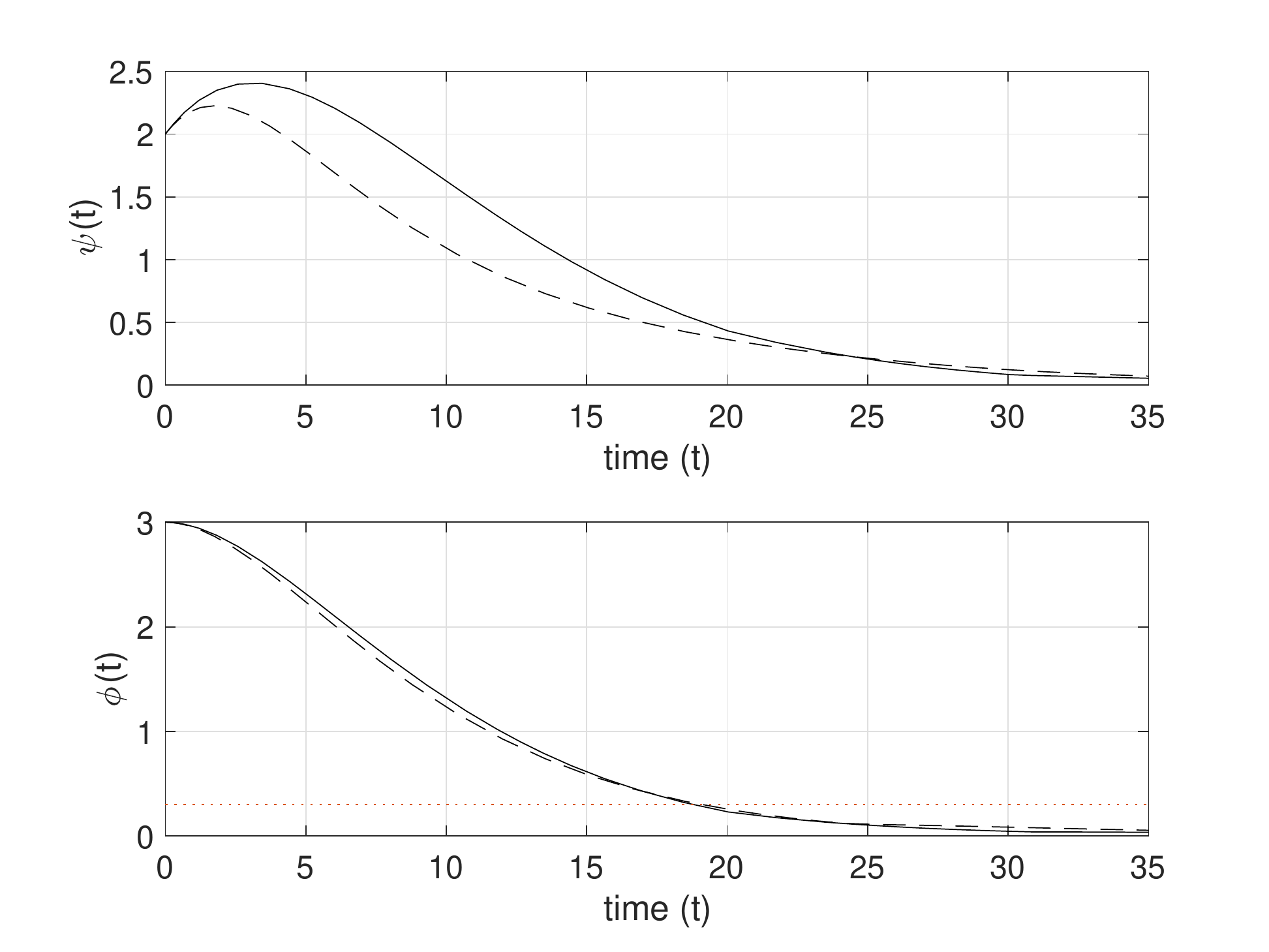}\vspace{-.3cm}
  \caption{Profile of agent difference $\psi(t)$ (top) and regulation error $\phi(t)$ (bottom); solid lines: $\rho =0.12$, $\varepsilon =0$;
dashed lines:  $\rho =0.2$, $\varepsilon =0.2$.}\label{fig4}
\end{figure}

Finally,  we show that the system performance can be modified by the dominant poles. Let $F =[ -0.2,\; 1]$. The corresponding results
are repeated in Figs.~\ref{fig3} and \ref{fig4}.  The aforementioned observation is still valid even when the closed-loop
system dynamics show more dampness.
\vspace{-.3cm}

\vspace{-.3cm}
\section{Conclusion}\label{s7}
 \vspace{-.1cm}
This paper has studied a consensus tracking problem of linear multi-agent systems. Compared with the conventional observer based control approach, where each observer estimates the leader's information via neighborhood communication, this paper features itself in two aspects: first, coupling effect between follower nodes in the physical layer  is deliberately introduced to take into account the cooperative behavior between all follower nodes  before they converge to the leader node; second, dominant poles of follower nodes can be adjusted to obtain a desired transient performance.
\vspace{-.2cm}
\section{Appendix}
Lemma \ref{F2} states the property of positive realness of a dynamic system under state feedback control.

\begin{lemma}[\cite{am}] \label{F2}
For the system
\begin{align*}
    \dot{x}=&Ax+Bu,
\end{align*}
suppose $(A,B)$ is stabilizable and the stabilizing feedback control gain is $F=B\t X$,
where $X>0$ is the stabilizing solution to the algebraic Riccati equation
\begin{align*}
    A\t X+XA-XBB\t X+{\mathcal Q}=0
\end{align*}
with ${\mathcal Q}> 0$.  Then the closed-loop transfer function
\begin{align*}
    T_F(s)=F(sI-A-BF)^{-1}B
\end{align*}
is positive real, i.e.,
$T_F(s)+T_{F}(s)^{\ast}\geq 0,\ \forall\ {\rm Re}[s]>0.$
\end{lemma}

 \vspace{-.1cm}
The result in Lemma \ref{FF} has been claimed in \cite{ygc2018}, and is summarized below with a self-contained proof.
\begin{lemma}  \label{FF}
If the matrices $M_1,M_2 \in \C^{N \times N}$ satisfy
\[M_1+M_1^*\geq 0,\quad M_2+M_2^* > 0,\]
then
$\det(I+M_1M_2)\neq0. $
\end{lemma}
 \vspace{-.2cm}
\bproof
Pick  $R_i =(I -M_i)(I + M_i)^{-1},\; i=1,2.$
One has
\begin{align*}
M_i=(I+R_i)^{-1}(I-R_i)=(I-R_i)(I+R_i)^{-1},\; i=1,2.
\end{align*}
Direct calculation shows
\begin{align*}
M_i+M_i^*
=&2(I+R_i)^{-1}(I-R_iR_i^*)(I+R_i^*)^{-1}.
\end{align*}
Then, $M_1+M_1^*\geq 0$ implies that the singular value $\overline{\sigma}(R_1)\leq 1$,
and  $M_2+M_2^* > 0$ implies $\overline{\sigma}(R_2) <  1$. Then, from
\begin{align*}
I+M_1M_2
=&2(I+R_1)^{-1}(I + R_1R_2)(I+R_2)^{-1},
\end{align*}
we have
\begin{align*}
\det(I+M_1M_2)=\frac{2 \det(I+R_1R_2)}{\det(I+R_1)\det(I+R_2)}\neq0,
\end{align*}
due to   $\overline{\sigma}(R_1R_2)\leq \overline{\sigma}(R_1)\overline{\sigma}(R_2)<1$.
\eproof

The following lemma is adopted from \cite{f}  with slight modification (cf. Theorem 4.1 and Proposition 4.1 of \cite{f}).

\begin{lemma}  \label{pole}
Consider the LQR problem \eqref{LQ} of the structure \eqref{AB2} under Assumption \ref{A1}.
Let  \begin{align}
 F = B_2^{-1} V^{-1} \Sigma V[K\;\; I] ,\; \Sigma = \mbox{diag} (\sigma_1, \cdots, \sigma_m),
\end{align}
where $\sigma_1, \cdots, \sigma_m>0$ and
$V, \Sigma, K$ are selected such that $w = -F z$ is a solution to the LQR problem.
In particular, the eigenvalues of $A_{11} -A_{12} K$ has the specified stable  eigenvalues
$\{s_{1}, \cdots, s_{n-m}\}$, none of which is an eigenvalue of $A_{11}$.
If $d \to \infty$,  then matrix $A-dBF$ has the following eigenvalue distribution:
\begin{itemize}

\item[(1)] there are $n-m$ eigenvalues $\lambda^j$ of the form $\lambda^j \to s_j$,
$j=1,\cdots, n-m$; and

\item[(2)] there are $m$ eigenvalues $\lambda^{j+n-m}$ of the form $\lambda^{j+n-m} \to -d \sigma_j$,
$j=1,\cdots, m$.

 \end{itemize}

 \end{lemma}

\vspace*{-0.25in}


\begin{thebibliography}{10}

\bibitem{am}
B.D.O. Anderson and J.B. Moore.
\newblock {\em Optimal Control -- Linear Quadratic Methods}.
\newblock Prentice-Hall, Englewoods Cliffs, NJ, 1990.


\bibitem{Rovithakis2017TAC-decentr}
C.P. Bechlioulis and G.A. Rovithakis.
\newblock {Decentralized robust synchronization of unknown high order nonlinear
  multi-agent systems with prescribed transient and steady state performance}.
\newblock {\em IEEE Transactions on Automatic Control}, 62(1):123--134, 2017.

\bibitem{CaiH2017Auto-adaptive}
H.~Cai, F.L. Lewis, G.~Hu, and J.~Huang.
\newblock {The adaptive distributed observer approach to the cooperative output
  regulation of linear multi-agent systems}.
\newblock {\em Automatica}, 75:299--305, 2017.

\bibitem{CaoYuRenChen2013Review}
Y.~Cao, W.~Yu, W.~Ren, and G.~Chen.
\newblock {An overview of recent progress in the study of distributed
  multi-agent coordination}.
\newblock {\em IEEE Transactions on Industrial Information}, 9(1):427--438,
  2013.

\bibitem{f}
T.~Fujii.
\newblock A new approach to the {LQ} design from the viewpoint of the inverse
  regulator problem.
\newblock {\em IEEE Transactions on Automatic Control}, 32(11):995--1004, 1987.

\bibitem{KnornChenMiddleton2016Overview}
S.~Knorn, Z.~Chen, and R.~Middleton.
\newblock {Overview: collective control of multi-agent systems}.
\newblock {\em IEEE Transactions on Control of Network Systems}, 3(4):334--347,
  2016.

\bibitem{Liang2019SMC-adaptive}
H.~Liang, Y.~Zhou, H.~Ma, and Q.~Zhou.
\newblock {Adaptive distributed observer approach for cooperative containment
  control of nonidentical networks}.
\newblock {\em IEEE Transactions on Cybernetics}, 49(2):299--307, 2019.

\bibitem{OhaParkAhn2015Survey}
K.K. Oha, M.C. Park, and H.S. Ahn.
\newblock {A survey of multi-agent formation control}.
\newblock {\em Automatica}, 53:424--440, 2015.

\bibitem{RenBook2008}
W.~Ren and R.W. Beard.
\newblock {\em {Distributed Consensus in Multi-Vehicle Cooperative Control:
  Theory and Applications}}.
\newblock Springer-Verlag, London, 2008.

\bibitem{Shim2015CDC-funnel}
H.~Shim and S.~Trenn.
\newblock {A preliminary result on synchronization of heterogeneous agents via
  funnel control}.
\newblock In {\em 54th IEEE Conference on Decision and Control {(CDC)}}, pages
  2229--2234, Dec. 2015.

\bibitem{SuYF2012TAC-output}
Y.~Su and J.~Huang.
\newblock Cooperative output regulation of linear multi-agent systems.
\newblock {\em IEEE Transactions on Automatic Control}, 57(4):1062--1066, 2012.

\bibitem{ygc2018}
F.~Yan, G.~Gu, and X.~Chen.
\newblock A new approach to cooperative output regulation for heterogeneous
  multi-agent systems.
\newblock {\em SIAM Journal on Control and Optimization}, 56(3):2074--2094,
  2018.

\bibitem{YuMaZhang2019SMC-prescribed}
T.~Yu, L.~Ma, and H.~Zhang.
\newblock {Prescribed performance for bipartite tracking control of nonlinear
  multi-agent systems with hysteresis input uncertainties}.
\newblock {\em IEEE Transactions on Cybernetics}, 49(4):1327--1338, 2019.

\bibitem{ZhangLewisDas2011TAC-optimalsyn}
H.~Zhang, F.L. Lewis, and A.~Das.
\newblock {Optimal design for synchronization of cooperative systems: state
  feedback, observer and output feedback}.
\newblock {\em IEEE Transactions on Automatic Control}, 56(8):1948--1952, 2011.

\bibitem{JL2020SMCA-State}
J.~Zhang, T.~Feng, H.~Zhang, and X.~Wang.
\newblock The decoupling cooperative control with dominant poles assignment.
\newblock {\em IEEE Transactions on Systems, Man, and Cybernetics: Systems}, to
  be published, DOI: 10.1109/TSMC.2020.3011142, 2020.

\end{thebibliography}

\end{document}